\documentclass[11pt,twoside,reqno]{amsart}

\usepackage{bm}

\usepackage{amssymb,amsmath,amstext,amsthm,amsfonts,amscd,xcolor}

\usepackage{dsfont}

\usepackage{amsthm, enumerate}

\usepackage[ansinew]{inputenc}
\usepackage{graphicx}
\usepackage[mathscr]{eucal}
\usepackage{hyperref}

\newcommand{\R}{\mathbb{R}}

\newcommand{\C}{\mathbb{C}}
\newcommand{\N}{\mathbb{N}}
\newcommand{\Z}{\mathbb{Z}}

\newcommand{\T}{\mathbb{T}}

\let\newpf\proof \let\proof\relax 
\newenvironment{pf}{\newpf[\proofname]}{\qed\endtrivlist}

\newtheorem{Theorem}{Theorem}[section]
\newtheorem*{Theorem*}{Theorem 1.1}
\newtheorem{Lemma}{Lemma}[section]
\newtheorem{Proposition}{Proposition}[section]

\newtheorem{Remark}{Remark}[section]

\newtheorem{Definition}{Definition}[section]

\newcommand{\normtwo}[1]{
	{\left\vert\kern-0.25ex\left\vert\kern-0.25ex\left\vert #1
		\right\vert\kern-0.25ex\right\vert\kern-0.25ex\right\vert} }

\newcommand{\la}{\langle}
\newcommand{\ra}{\rangle}

\makeatletter
\newsavebox\myboxA
\newsavebox\myboxB
\newlength\mylenA

\newcommand*\xoverline[2][0.75]{%
	\sbox{\myboxA}{$\m@th#2$}%
	\setbox\myboxB\null
	\ht\myboxB=\ht\myboxA%
	\dp\myboxB=\dp\myboxA%
	\wd\myboxB=#1\wd\myboxA
	\sbox\myboxB{$\m@th\overline{\copy\myboxB}$}
	\setlength\mylenA{\the\wd\myboxA}
	\addtolength\mylenA{-\the\wd\myboxB}%
	\ifdim\wd\myboxB<\wd\myboxA%
	\rlap{\hskip 0.5\mylenA\usebox\myboxB}{\usebox\myboxA}%
	\else
	\hskip -0.5\mylenA\rlap{\usebox\myboxA}{\hskip 0.5\mylenA\usebox\myboxB}%
	\fi}
\makeatother

\newcommand\restr[2]{{
		\left.\kern-\nulldelimiterspace 
		#1 
		\vphantom{\big|} 
		\right|_{#2} 
}}

\newcommand{\diam}{\mathrm{diam}}

\begin{document}
	
	\title{Quantitative Reducibility of $\bf{C^k}$ Quasi-Periodic Cocycles}
	
	\author{Ao Cai, Huihui Lv, Zhiguo Wang}
		 
	\begin{abstract}
	This paper establishes an extreme $C^k$ reducibility theorem of quasi-periodic $SL(2, \R)$ cocycles in the local perturbative region, revealing both the essence of Eliasson [Commun.Math.Phys.1992] and Hou-You [Invent.Math.2012] in respectively the non-resonant and resonant cases. By paralleling further the reducibility process with the almost reducibility, we are able to acquire the least initial regularity as well as the least loss of regularity for the whole KAM iterations. This, in return, makes various spectral applications of quasi-periodic Schr\"odinger operators wide open.	
	\end{abstract}
	
	\maketitle

	\section{Introduction}
	Assume that $(M, \mathcal{B}, \mu)$ is a probability space and $f: M \to M$ is an invertible map which preserves the measure $\mu$ and is ergodic with respect to $\mu$. Let $A: M \to SL(2, \R)$ be a measurable function. The linear cocycle defined by $A$ over the base dynamics $f$ is the transformation: 
	$$
	(f,A):M \times \R^2 \to M \times \R^2; \ (\theta, v) \mapsto (f(\theta), A(\theta) \cdot v). 
	$$
	Note that the iterates of $(f,A)$ have the form $(f,A)^n = (f^n, A_n)$, where $A_n(\theta) = A(f^{n-1}(\theta)) \cdots A(f(\theta)) A(\theta), n \in \N$ and $A_{-n}(\theta) = A_n(f^{-n}(\theta))^{-1}$. In particular, a $C^k$ quasi-periodic linear cocycle $(\alpha, A)$ consists of a rationally independent $\alpha \in \T^d$, which determines an ergodic torus translation on the base and $A \in C^k(\T^d, SL(2,\R))$ which is a $k$ times differentiable matrix-valued function with continous $k$-th derivatives. 
	
	In this paper, we present a purified quantitative reducibility theorem for finitely differentiable quasi-periodic cocycles. In particular, it applies to the $C^k$ quasi-periodic Schr\"{o}dinger cocycles $(\alpha, A)$, where
	$$
	A(\theta)=S_E^{V}(\theta)=\begin{pmatrix} E- V(\theta) & -1 \\ 1 & 0 \end{pmatrix}. 
	$$
	This effectively extends the spectral applications in the one-dimensional discrete quasi-periodic Schr\"{o}dinger operator with a $C^k$ potential: 
	\begin{equation}\label{so}
		(H_{V,\alpha,\theta}x)_n=x_{n+1}+x_{n-1}+V(\theta+n\alpha)x_{n}, n\in \Z,
	\end{equation}
	as any formal solution (not necessarily in $\ell^2$) of $H_{V,\alpha,\theta}x=Ex$ satisfies 
	$$
	A(\theta+n\alpha)\begin{pmatrix} x_{n} \\ x_{n-1} \end{pmatrix}=\begin{pmatrix} x_{n+1} \\ x_{n} \end{pmatrix}.
	$$
	In equation (1), $V\in C^k(\T^d,\R)$ is called the potential, $\alpha \in \R^d$ is called the frequency satisfying $\la n,\alpha\ra \notin \Z$ for any $n\in \Z^d \backslash\{0\}$ and $\theta\in \T^d=\R^d/\Z^d$ (or $\R^d/(2\pi\Z)^d$ if preferable) is called the initial phase, $k,d\in \N^+$.
	The spectrum of $H_{V,\alpha,\theta}$ is denoted as $\Sigma_{V,\alpha,\theta}$, which is independent of the phase $\theta$. When $V$ is bounded, $H_{V,\alpha,\theta}$ is a bounded self-adjoint operator on $\ell^2(\Z)$ and $\Sigma_{V,\alpha,\theta} \subset \R$ is a compact perfect set. See the nice survey of Damanik \cite{MR3681983} for more details of the Schr\"odinger operator unity. Moreover, readers are invited to consult the excellent book of Damanik and Fillman \cite{damanik2022one} which is significant and timely for the community.

	\subsection{Quantitative reducibility}
	The reducibility of quasi-periodic cocycles aims to conjugate the original quasi-periodic cocycles to constant ones via transformations that are essentially coordinate changes. However, due to topological obstructions, reducibility may fail and the concept of almost reducibility naturally arises, which settles for almost conjugating to constant cocycles but has proven to be very powerful in studying the spectral theory of quasi-periodic Schr\"odinger operators \cite{you2018quantitative}. Readers are invited to Section 2 for precise definitions of (almost) reducibility. 
	
	In the previous literature, almost reducibility is regarded as a prerequisite for reducibility, especially in the analytic region, due to the great flexibility of shrinking the analytic radius arbitrarily. Nevertheless, in the $C^k$ topology, the least loss of regularity is a fixed number. In order to establish a $C^k$ reducibility theorem with least initial regularity and most remainder, we perform a process parallel to the process of almost reducibility. In other words, reducibility is built not after almost reducibility, but at the same time with extra assumptions of the fibered rotation number $\rho(\alpha, A)$.  
	
	Recall that $\alpha \in\R^d$ is called {\it Diophantine} if there exist $\kappa>0$ and $\tau>d$ such that $\alpha \in {\rm DC}_d(\kappa,\tau)$, where
	\begin{equation}\label{dio1}
		{\rm DC}_d(\kappa,\tau):=\left\{\alpha \in\R^d:  \inf_{j \in \Z}\left|\la m,\alpha  \ra - j \right|
		> \frac{\kappa}{|m|^{\tau}},\quad \forall \, m\in\Z^d\backslash\{0\} \right\}.
	\end{equation}
	Here we denote
	$$
	\lvert m\rvert=\lvert m_1\rvert+\lvert m_2\rvert+ \cdots + \lvert m_d\rvert
	$$
	and
	$$
	\la m,\alpha \ra =m_1\alpha_1+m_2\alpha_2+\cdots+m_d\alpha_d.
	$$
	Denote ${\rm DC}_d=\bigcup_{\kappa,\tau} {\rm DC}_d(\kappa,\tau)$, which is of full Lebesgue measure.
	
	The rotation number $\rho (\alpha, A)$ is called {\it Diophantine} w.r.t. $\alpha$ if it satisfies the condition $\rho (\alpha, A) \in {\rm DC}^{\alpha}_d(\gamma,\tau)$, where $\gamma > 0$, $\tau > d$ and
	\begin{equation}\label{dio2}
		{\rm DC}^{\alpha}_d(\gamma,\tau):=\left\{\phi \in\R:  \inf_{j \in \Z}\left|2\phi - \la m,\alpha  \ra - j \right|
		> \frac{\gamma}{(|m|+1)^{\tau}},\quad \forall \, m\in\Z^d \right\}.
	\end{equation}
	The rotation number $\rho (\alpha, A)$ is called rational w.r.t. $\alpha$ if $2\rho (\alpha,A) = \la m_0,\alpha \ra$ mod $\Z$ for some $m_0 \in \Z^d$.
	
	Our main theorem is as follows: 
	
	\begin{Theorem}\label{main1}
		Let $A \in SL(2, \R), \alpha \in \rm DC_d (\kappa, \tau) $ and $f \in C^k(\T^d, sl(2, \R))$ with $k > 14\tau + 2$. There exists $\epsilon = \epsilon(\kappa, \tau, d, k, ||A||)$ such that if $\lVert f(\theta) \rVert_k \leqslant \epsilon$
		and $\rho(\alpha, Ae^f)$ is Diophantine or rational with respect to $\alpha$, then $(\alpha, Ae^{f(\theta)})$ is $C^{k,k_0}$ reducible with $k_0 < k-10\tau-3$.
	\end{Theorem}
	
	\begin{Remark}
		The quantitative version of Theorem 1.1 can be found in Section 3, specifically in Theorem \ref{thm3.1}. Our theorem is a significant improvement over the results presented in \cite{cai2022reducibility}. Notably, the loss of regularity $10\tau+3$ does not depend on the parameter $k$.
	\end{Remark}
	
	As there are few reducibility results in the $C^k$ topology, we only briefly review the development of analytic reducibility. Let's first discuss it for local cocycles of the form $(\alpha, Ae^{f(\theta)})$. An achievement was initiated by Dinaburg and Sinai \cite{dinaburg1975one}, who applied the classical KAM scheme. They established the positive measure reducibility in terms of the rotation number for continuous quasi-periodic Schr\"{o}dinger equations featuring small analytic potentials in the perturbative regime. Moser and P\"{o}schel \cite{moser1984extension} later expanded on this achievement by utilizing a resonance-cancellation technique, which extended the positive measure reducibility to a class of rotation numbers that are rational w.r.t. $\alpha$. Moreover, the breakthrough came from Eliasson \cite{MR1167299}, who established weak almost reducibility for all energies $E$ and full measure reducibility for Diophantine frequencies and small analytic potentials. For further insights into strong almost reducibility results, readers can refer to the works of Chavaudret \cite{chavaudret2013strong} and Leguil-You-Zhao-Zhou \cite{leguil2017asymptotics}. In the non-perturbative regime, Puig \cite{puig2005nonperturbative} employed the localization method to derive a non-perturbative version of Eliasson's reducibility theorem. As for continuous linear systems, Hou-You \cite{hou2012almost} proved weak almost reducibility results for all rotation numbers and frequencies $\omega = (1, \alpha) \in \mathbb{T}^2$ under small analytic perturbations.
	
	In the case of global cocycle $(\alpha, A(\theta))$ with $A \in C^\omega(\T, SL(2, \R))$ and $\alpha \in \T$, Avila and Krikorian \cite{2006Reducibility} applied the renormalization scheme to show that, for $\alpha$ satisfying certain recurrent Diophantine conditions and almost every $E$, the quasi-periodic Schr\"{o}dinger cocycle is either reducible or non-uniform hyperbolic. Additionally, Avila, Fayad and Krikorian \cite{avila2011kam} proved that for irrational $\alpha$ and almost every $E$, the quasi-periodic Schr\"{o}dinger cocycle is either rotations reducible or non-uniformly hyperbolic. 
	
	In the realm of finitely differentiable topology, reducibility results are scarce and the existing ones are rough in the sense of very large initial regularity and very much loss. We note that our philosophy is that nicer reducibility implies nicer spectral applications. In this direction, Theorem 1.1, which was established via remaining all the parameters, allows us to obtain the most general applications according to the technique.

	\subsection{Spectral type and structure}
	Over the past forty years, advances have been made in studying the spectral theory of the Schr\"{o}dinger operator, focusing on understanding the spectral type and the structure of the spectrum. In the 21st century, people found that quantitative dynamical estimates lead to quantitative spectral applications \cite{MR2578605}. In particular, quantitative (almost) reducibility is one of the most powerful techniques, with abundant fruitful spectral results \cite{avila2009cantor}\cite{cai2022absolutely}\cite{cai2019sharp}\cite{cai2022reducibility}\cite{cai2021polynomial}\cite{damanik2016spectrum}\cite{jian2019sharp}\cite{johnson1982rotation}\cite{puig2005nonperturbative}\cite{wang2017cantor}. One is invited to consult You’s 2018 ICM survey \cite{you2018quantitative}. Generally, spectral type results refer to absolutely continuous spectrum, singular continuous spectrum, pure point spectrum, ballistic transport, Anderson localization, etc. while regarding spectral structure, common results often involve the Cantor spectrum and homogeneous spectrum. For the sake of concision, we only list one spectral type application and one spectral structure application respectively, though our Theorem 1.1 along with its quantitative version is promisingly applicable to many other spectral results in relation. 
	
	\begin{Theorem}(Spectral type application)\label{main2}
		Assume $\alpha\in {\rm DC}_d(\kappa, \tau)$, $V\in C^k(\T^d, \R)$ with $k>14\tau+2$. If there exists $\epsilon'=\epsilon'(\kappa, \tau, k, d)$ such that $\lVert V \rVert_k \leqslant \epsilon'$ and $\rho(\alpha, S_E^V)$ is Diophantine or rational w.r.t. $\alpha$ for $E \in \Sigma_{V,\alpha}$, then $H_{V,\alpha,\theta}$ has strong ballistic transport for all $\theta \in \T^d$. 
	\end{Theorem}
	
	\begin{Remark}
		The definition of (global) strong ballistic transport is too lengthy, so we prefer to put it in Section 4 rather than here. This theorem compensates for the work of Ge and Kachkovskiy \cite{ge2023ballistic} by providing the precise requirement of $k$ for a $C^k$ quasi-periodic Schr\"{o}dinger family to be reducible.  
	\end{Remark}
	
	Indeed, there is a qualitative connection between spectral type and transport properties, which can be understood through the RAGE theorem \cite{MR3681983}. However quantitative results also illustrate the connection between them. For example, except for purely absolutely continuous spectra, we generally don't expect to observe ballistic transport phenomena. In particular, it has been proven in \cite{simon1990absence} that point spectra cannot support any ballistic motion. However, when the operator is restricted to a subspace that supports purely absolutely continuous spectrum, we can still expect to observe ballistic transport phenomena. And the works of Guarneri, Combes and Last \cite{combes1993connections} \cite{last1996quantum} have also made contributions to the connection. In particular, the Guarneri-Combes-Last theorem \cite{last1996quantum} quantitatively provides insights into transport phenomena in one-dimensional systems with absolutely continuous spectrum. In subsequent studies, Asch and Knauf \cite{asch1997motion}, as well as Damanik \cite{damanik2015quantum}, demonstrated the occurrence of strong ballistic transport in periodic continuous Schrödinger operators, which are widely known to have absolutely continuous spectra. Later, Zhang and Zhao \cite{Zhao2017Ballistic} explicitly established a connection between the values of transport exponents and absolutely continuous spectra in the setting of discrete single-frequency quasi-periodic operators. More recently, based on the results and techniques in the discrete setting of Fillman \cite{fillman2017ballistic},  Giorgio's research \cite{2021Ballistic} discovered strong ballistic transport in a class of continuum limit-periodic operators known to possess absolutely continuous spectra. 
	
	Let us now state our spectral structure result. First, recall that
	
	\begin{Definition}\cite{cai2021polynomial}
		Let $\nu > 0$, a closed set $\mathfrak{B} \subset \R$ is called $\nu$-homogeneous if 
		$$
		|\mathfrak{B} \cap (E-\epsilon, E+\epsilon)|>\nu \epsilon, \ \forall E\in \mathfrak{B}, \ \forall 0<\epsilon<\diam \mathfrak{B}. 
		$$
	\end{Definition}
	
	As another corollary of Theorem \ref{main1}, we have 
	
	\begin{Theorem}(Spectral structure application)\label{main3}
		Assume $\alpha\in {\rm DC}_d(\kappa, \tau)$, $V\in C^k(\T^d, \R)$ with $k>17\tau+2$. If there exists $\epsilon''=\epsilon''(\kappa, \tau, k, d)$ such that $\lVert V \rVert_k \leqslant \epsilon''$, then $\Sigma_{V,\alpha}$ is $\nu$-homogeneous for some $\nu \in (0,1)$. 
	\end{Theorem}
	
	\begin{Remark}
		We know the homogeneity of the spectrum is related to polynomial decay of gap length and H\"{o}lder continuity of the integrated density of states. Therefore, we need to change the assumption to $k > 17\tau + 2$ so that we can further ensure the $\frac{1}{2}$-H\"{o}lder continuity of IDS (see Theorem 3.3 and Theorem 3.4 in \cite{cai2022absolutely}). This theorem greatly reduces the initial regularity assumption in the $C^k$ case.  
	\end{Remark}
    
    The spectrum's homogeneity plays a crucial role in the inverse spectral theory, as demonstrated in the seminal works of Sodin and Yuditskii \cite{2013Exponential}\cite{1996Dimensional}.  Under the assumption of a finite total gap length and a reflectionless condition on the spectrum, it has been proved that the homogeneity of the spectrum implies the almost periodicity of the associated potentials \cite{2002Absolutely}. In particular, the homogeneity of the spectrum is closely linked to Deift's conjecture, which investigates whether the solutions of the KdV equation exhibit quasi-periodicity when the initial data is quasi-periodic \cite{MR2578605}\cite{2002On}. In the continuous setting, Binder-Damanik-Goldstein-Lukic \cite{MR2578605} demonstrated that when considering small analytic quasi-periodic initial data with Diophantine frequency, the solution of the KdV equation exhibits almost periodicity in the temporal variable. In the discrete setting, Leguil-You-Zhao-Zhou \cite{2019Exponential} demonstrated that for the subcritical potential $V \in C^\omega(\T, \R)$, the Toda flow is almost periodic in the time variable when considering initial data that are also almost periodic with $\beta(\alpha) = 0$. Lately, Avila-Last-Shamis-Zhou \cite{avila2015global} also constructed an intriguing counter-example that even for the AMO, its spectrum is not homogeneous if $e^{-\frac{2}{3}\beta(\alpha)} < \lambda < e^{\frac{2}{3}\beta(\alpha)}$. In the $C^k$ case, Cai and Wang \cite{cai2021polynomial} have recently proved the homogeneity of the spectrum. They achieved this through a not so refined quantitative $C^k$ reducibility theorem for quasi-periodic $SL(2, \R)$ cocycles, as well as by employing the Moser-P\"{o}schel argument for the associated Schrödinger cocycles.

	\subsection{Structure of the paper}
	We formulate the structure of our paper as follows.
	In Section 2, we give some useful definitions and notations. In Section 3, we prove the quantitative reducibility of $C^k$ quasi-periodic $SL(2,\R)$ cocycles with the rotation number being Diophantine and rational w.r.t. the frequency in the local perturbative regime. As spectral applications, in Section 4, we use the $C^k$ quantitative reducibility theorem to get some results about ballistic transport and homogeneous spectrum in our settings.

	\section{Preliminaries}

	\subsection{Conjugation and reducibility}
	For a bounded analytic (possibly matrix valued) function $F(\theta)$ defined on $\mathcal{S}_h:= \{ \theta=(\theta_1,\dots, \theta_d)\in \mathbb{C}^d\ |\ \ | \Im \theta_i | < h, \forall 1\leqslant i\leqslant d\}$, let $|F|_h=  \sup_{\theta\in \mathcal{S}_h } \| F(\theta)\| $ and denote by $C^\omega_{h}(\T^d,*)$ the set of all these $*$-valued functions ($*$ will usually denote $\R$, $sl(2,\R)$, $SL(2,\R)$). Let $C^\omega(\T^d,*)=\cup_{h>0}C^\omega_{h}(\T^d,*)$ and $C^{k}(\T^{d},*)$ be the space of $k$ times differentiable with continuous $k$-th derivatives functions. Define the norm as
	$$
	\lVert F \rVert _{k}=\sup_{\substack{
			\lvert k^{'}\rvert\leqslant k,
			\theta \in \T^{d}
	}}\lVert \partial^{k^{'}}F(\theta) \rVert.
	$$ 
	
	For two cocycles $(\alpha,A_1)$, $(\alpha,A_2)\in \T^d  \times C^{\ast}(\T^d,SL(2,\R))$, ``$\ast$'' represents ``$\omega$'' or ``$k$'',  we can say that they are $C^{\ast}$ conjugated if there exists $Z\in C^{\ast}(2\T^d, SL(2,\R))$, such that $$
	Z(\theta+\alpha)A_1(\theta)Z^{-1}(\theta)=A_2(\theta).
	$$
	Notably, we want to define $Z$ on the $2\T^d=\R^d/(2\Z)^d$ for the purpose of making it still real-valued.
	
	An analytic cocycle $(\alpha, A)\in \T^d    \times C^{\omega}_h(\T^d, SL(2,\R))$ is called almost reducible if there exist a sequence of constant matrices $A_j\in SL(2,\R)$, a sequence of conjugations $Z_j\in C^{\omega}_{h_j}(2\T^d, SL(2,\R))$ and a sequence of small perturbation $f_j \in C^{\omega}_{h_j}(\T^d, sl(2,\R))$ such that
	$$
	Z_j(\theta+\alpha)A(\theta)Z_j^{-1}(\theta)=A_j e^{f_j(\theta)}
	$$
	with
	$$
	\lvert f_j(\theta)\rvert_{h_j}\rightarrow 0, \ \ j\rightarrow \infty.
	$$
	Furthermore, it is said to be weak $(C^{\omega})$ almost reducible if $h_j \rightarrow 0$ and it is said to be strong $(C^{\omega}_{h_j,h'})$ almost reducible if $h_j\rightarrow h'>0$. We also claim $(\alpha, A)$ is $C^{\omega}_{h,h'}$ reducible if there exist a constant matrix $\widetilde{A} \in SL(2,\R)$ and a conjugation map $\widetilde{Z}\in C^{\omega}_{h'}(2\T^d, $ $SL(2,\R))$ such that
	$$
	\widetilde{Z}(\theta+\alpha)A(\theta)\widetilde{Z}^{-1}(\theta)=\widetilde{A}(\theta).
	$$
	
	To avoid repetition, we have provided an equivalent definition of $C^k$ (almost) reducibility as follows. 
	
	One can say that a finitely differentiable cocycle $(\alpha,A)$ is $C^{k,k_1}$ almost reducible, if $A\in C^k(\T^d,SL(2,\R))$ and the $C^{k_1}$-closure of its $C^{k_1}$ conjugacies contains a constant. Besides, we say $(\alpha,A)$ is  $C^{k,k_1}$ reducible, if $A\in C^k(\T^d,SL(2,\R))$ and its $C^{k_1}$ conjugacies contain a constant.

	\subsection{Rotation number and degree}
	Suppose $A\in C^0(\T^d,SL(2,\R))$ is homotopic to identity. Then we show the projective skew-product $F_A:\T^d \times \mathbb{S}^1 \rightarrow \T^d \times \mathbb{S}^1$ with
	$$
	F_A(x,\omega):=\left(   x+\alpha, \frac{A(x)\cdot \omega}{\lvert A(x)\cdot \omega\rvert}\right),
	$$
	which is homotopic to identity as well. Thus, we will lift $F_A$ to a map $\widetilde{F}_A:\T^d\times \R\rightarrow \T^d\times \R$ with  $\widetilde{F}_A(x,y)=(x+\alpha,y+\psi(x,y))$, where for every $x \in \T^d$, $\psi(x,y)$ is $\Z$-periodic in $y$. The map $\psi:\T^d \times \R \rightarrow \R$ is said to be a lift of $A$. Assume $\mu$, which is invariant by $\widetilde{F}_A$, be any probability measure on $\T^d\times \R$. Its projection on the first coordinate is provided by the Lebesgue measure. The number
	\begin{equation}\label{rot1}
		\rho_{(\alpha,A)}:=\int_{\T^d\times \R} \psi(x,y)d\mu(x,y)\mbox{mod}\,\Z
	\end{equation}
	has nothing to do with the choices of the lift $\psi$ or the measure $\mu$. One calls it the {\it fibered rotation number} of cocycle $(\alpha, A)$ (readers can refer to \cite{johnson1982rotation} for more details).
	
	Assume
	$$
	R_{\phi}:=\begin{pmatrix} \cos 2\pi \phi & -\sin 2\pi \phi \\ \sin 2\pi \phi & \cos 2\pi \phi \end{pmatrix},
	$$
	if $A\in C^0(\T^d,SL(2,\R))$ is homotopic to $\theta\rightarrow R_{\la n,\theta\ra}$ for some $n\in \Z^d$, then we call $n$ the {\it degree} of $A$ and denote it by deg$A$. Furthermore, 
	\begin{equation}\label{deg1}
		\deg(AB)=\deg A+\deg B.
	\end{equation}
	
	Please note that the fibered rotation number remains invariant under real conjugacies that are homotopic to the identity map. Generally speaking, when the cocycle $(\alpha,A_1)$ is conjugated to $(\alpha,A_2)$ by $B\in C^0(2\T^d, SL(2,\R))$, i.e. $B(\cdot +\alpha)A_1(\cdot)B^{-1}(\cdot)=A_2(\cdot)$, we have
	\begin{equation}\label{rot2}
		\rho_{(\alpha,A_2)}=\rho_{(\alpha,A_1)} - \frac{\la \deg B,\alpha\ra}{2}.
	\end{equation}

	\subsection{Hyperbolicity and integrated density of states}
	We call the cocycle $(\alpha,A)$ {\it uniformly hyperbolic} if for every $\theta \in \T^d$, there exists a continuous decomposition $\C^2 = E^s(\theta) \oplus E^u(\theta)$ such that for some constants $C > 0$, $ c > 0$ and every $n \geqslant 0$,
	$$
	|A_n(\theta)v| \leqslant Ce^{-cn}|v|, v \in E^s(\theta),
	$$
	$$
	|A_{-n}(\theta)v| \leqslant Ce^{-cn}|v|, v \in E^u(\theta).
	$$
	This decomposition is invariant by the dynamics, which means that for any $\theta \in \T^d$, $A(\theta)E^*(\theta) = E^*(\theta + \alpha)$, for $*= s, u$. In the $C^0$ topology, the set of uniformly hyperbolic cocycles is an open set. Specifically, in the case of quasi-periodic Schr\"{o}dinger operators, the cocycle $(\alpha, S_E^V)$ is uniformly hyperbolic if and only if $E \notin \Sigma_{V,\alpha}$, or in other words, if the energy lies within a spectral gap \cite{johnson1986exponential}.
	
	Let's consider the Schr\"{o}dinger operators $H_{V,\alpha,\theta}$, where an important concept is the integrated density of states (IDS). The IDS is a function $N_{V,\alpha}: \mathbb{R}\rightarrow [0,1]$ that can be defined by 
	\[
	N_{V,\alpha}(E) = \int_{\mathbb{T}^{d}}\mu_{V,\alpha,\theta}(-\infty,E]d\theta,
	\]
	where  $\mu_{V,\alpha,\theta} = \mu_{V,\alpha,\theta}^{e_{-1}}+\mu_{V,\alpha,\theta}^{e_{0}}$ is the universal spectral measure of $H_{V,\alpha,\theta}$ and  $\{e_{i}\}_{i\in \Z}$ is the cannonical basis of $\ell^{2}(\mathbb{Z})$. We say $\{e_{-1}, e_{0}\}$ is the pair of cyclic vectors of $H_{V,\alpha,\theta}$ here.
	
	There exist alternative approaches to defining the IDS by counting eigenvalues of the truncated Schr\"{o}dinger operator. For further details, readers may consult \cite{avron1983almost}. Furthermore, there is a connection between $\rho(\alpha, S^V_E)$ and the IDS, which can be expressed as follows:
	\begin{equation}\label{ids}
		N_{V,\alpha}(E)=1-2\rho(\alpha,S^{V}_{E})  \mod \mathbb{Z}. 
	\end{equation}
	
	To gain a better understanding of the various types of spectral measures, I recommend referring to the book \cite{MR3681983}.

	\subsection{Analytic approximation}
	Let $f \in C^{k}(\T^{d},sl(2,\R))$. By Zehnder \cite{zehnder1975generalized}, there is a sequence $\{f_{j}\}_{j\geqslant 1}$, $f_{j}\in C_{\frac{1}{j}}^{\omega}(\T^{d},sl(2,\R))$ and a universal constant $C'$, such that
	\begin{eqnarray} \nonumber \lVert f_{j}-f \rVert_{k} &\rightarrow& 0 , \quad  j \rightarrow +\infty, \\
		\label{aa}\lvert f_{j}\rvert_{\frac{1}{j}} &\leqslant& C'\lVert f \rVert_{k}, \\   \nonumber \lvert f_{j+1}-f_{j} \rvert_{\frac{1}{j+1}} &\leqslant& C'(\frac{1}{j})^k\lVert f \rVert_{k}.
	\end{eqnarray}
	Furthermore, if $k\leqslant \tilde{k}$ and $f\in C^{\tilde{k}}$, the properties $(\ref{aa})$ still hold with $\tilde{k}$ instead of $k$. This implies that the sequence can be constructed from $f$ irrespective of its regularity (since $f_{j}$ is achieved by convolving $f$ with a map that does not depend on $k$).

	\section{Dynamical estimates: full measure reducibility}\label{sec3}
	In this section, the main emphasis is on investigating the reducibility property of the following $C^k$ quasi-periodic $SL(2,\R)$ cocycle: 
	$$
	(\alpha,Ae^{f(\theta)}): \T^d \times \R^2 \to \T^d \times \R^2; (\theta, v) \mapsto (\theta + \alpha, Ae^{f(\theta)} \cdot v).
	$$
	where $\alpha\in {\rm DC}_d(\kappa,\tau)$, $A\in SL(2,\R), \, f\in C^k(\T^{d},sl(2,\R))$ and $d\in \N^+$. Our approach involves initially analyzing the approximating analytic cocycles $\{(\alpha,Ae^{f_{j}(\theta)})\}_{j\geqslant 1}$ and subsequently transferring the obtained estimates to the targeted $C^k$ cocycle $(\alpha,Ae^{f(\theta)})$ through analytic approximation techniques.

	\subsection{Preparations}
	In the following subsections, we will consider fixed parameters $\rho,\epsilon,N,\sigma$; we refer to the situation in which there exists $m_\ast$ satisfying $0<\lvert m_\ast\rvert \leqslant N$ such that
	$$
	\inf_{j \in \Z}\lvert 2\rho- \la m_\ast,\alpha\ra - j\rvert< \epsilon^{\sigma},
	$$
	as the ``resonant case'' (for simplicity, we will use the notation ``$\lvert 2\rho - \la m_\ast,\alpha\ra \rvert$" to represent the left side and similarly, ``$\lvert \la m_\ast,\alpha\ra \rvert$" to represent the right side). The integer vector $m_\ast$ will be referred to as a ``resonant site." This type of small divisor problem commonly arises when attempting to solve the cohomological equation at each step of the KAM procedure. 
	Resonances are connected to a useful decomposition of the space $\mathcal{B}_r:=C^{\omega}_{r}(\T^{d},su(1,1))$. For further details and the precise definition of $su(1,1)$ and $SU(1,1)$, please refer to the proof of Proposition $\ref{pro3.1}$. 
	
	Given $\alpha\in \R^{d}$, $A\in SU(1,1)$ and $\eta>0$, there exists a decomposition $\mathcal{B}_r=\mathcal{B}_r^{nre}(\eta) \bigoplus\mathcal{B}_r^{re}(\eta)$ satisfying that for any $Y\in\mathcal{B}_r^{nre}(\eta)$,
	\begin{equation}\label{space}
		A^{-1}Y(\theta+\alpha)A\in\mathcal{B}_r^{nre}(\eta),\,\lvert A^{-1}Y(\theta+\alpha)A-Y(\theta)\rvert_r\geqslant\eta\lvert Y(\theta)\rvert_r.
	\end{equation}
	And denote $\mathbb{P}_{nre}$, $\mathbb{P}_{re}$ the standard projections from $\mathcal{B}_r$ onto $\mathcal{B}_r^{nre}(\eta)$ and $\mathcal{B}_r^{re}(\eta)$ respectively.
	
	Next, we have a crucial lemma that plays a key role in eliminating all the non-resonant terms. This lemma will be utilized in the proof of the resonant conditions.
	
	\begin{Lemma}\label{lem3.1}\cite{cai2019sharp}\cite{hou2012almost}
		Assume that $A\in SU(1,1)$, $\epsilon\leqslant (4\lVert A\rVert)^{-4}$ and  $\eta \geqslant 13\lVert A\rVert^2{\epsilon}^{\frac{1}{2}}$. For any $g\in \mathcal{B}_r$ with $|g|_r \leqslant \epsilon$,  there exist $Y\in \mathcal{B}_r$ and $g^{re}\in \mathcal{B}_r^{re}(\eta)$ such that
		$$
		e^{Y(\theta+\alpha)}(Ae^{g(\theta)})e^{-Y(\theta)}=Ae^{g^{re}(\theta)},
		$$
		with estimates
		$$\lvert Y \rvert_r\leqslant \epsilon^{\frac{1}{2}}, \ \lvert g^{re}\rvert_r\leqslant 2\epsilon.
		$$
	\end{Lemma}
	
	\begin{Remark}\label{rem3.1}
		For the inequality ``$\eta \geqslant 13\lVert A\rVert^2{\epsilon}^{\frac{1}{2}}$'', ``$\frac{1}{2}$'' is sharp because of the quantitative Implicit Function Theorem \cite{berti2006forced}\cite{deimling1989nonlinear}. The proof relies solely on the fact that $\mathcal{B}_r$ is a Banach space. Therefore, it is applicable not only to the $C^{\omega}$ topology but also to the $C^k$ and $C^0$ topologies. For more detailed information, please refer to the appendix of \cite{cai2019sharp}. 
	\end{Remark}

	\subsection{Analytic KAM theorem}
	Following our plan, our first objective is to establish the KAM theorem for the analytic quasi-periodic $SL(2,\R)$ cocycle $(\alpha, Ae^{f(\theta)})$, where $A$ possesses eigenvalues ${e^{i\rho},e^{-i\rho}}$ with $\rho \in 2\pi\R \cup 2\pi i\R$. We present our quantitative analytic KAM theorem in the following. 
	 
	\begin{Proposition}\label{pro3.1}\cite{cai2022absolutely}\cite{cai2019sharp}
		Let $\alpha \in {\rm DC}_d(\kappa,\tau), \kappa,r>0, \tau>d,\sigma<\frac{1}{6}$. Suppose that $A \in SL(2,\R)$ satisfying  $||A||$ bounded, $f \in C_r^{\omega}(\T^{d},sl(2,\R))$. Then for any $0<r'<r$, there exist constants $c=c(\kappa,\tau,d), D>\frac{2}{\sigma}$ and $\tilde{D}=\tilde{D}(\sigma)$ such that if
		\begin{equation}\label{estrf}
			\lvert f \rvert_r\leqslant\epsilon \leqslant \frac{c}{\lVert A\rVert^{\tilde{D}}}(r-r')^{D\tau},
		\end{equation}
		then there exist $B \in C^{\omega}_{r'}(2\T^d,SL(2,\R)), A_+ \in SL(2,\R)$ and $f_+ \in C^{\omega}_{r'}(\T^d, \\ sl(2,\R))$ such that
		$$
		B(\theta+\alpha)(Ae^{f(\theta)})B^{-1}(\theta)=A_+e^{f_+(\theta)}.
		$$
		More precisely, let $N= \frac{2}{r-r'} \lvert \ln \epsilon \rvert$,then we can distinguish two cases: 
		\begin{itemize}
			\item (Non-resonant case) if for any $m \in \Z^d$ with $0<|m| \leqslant N$, we have 
			$$
			\lvert 2 \rho-\langle m,\alpha \rangle \rvert \geqslant {\epsilon}^{\sigma},
			$$
			then
			$$
			\lvert B(\theta)- Id\rvert_{r'}\leqslant \epsilon^{1-\frac{8}{D}}, \lvert f_+\rvert_{r'}\leqslant \epsilon^{2-\frac{8}{D}}
			$$
			and
			$$
			\lVert A_+ - A\rVert \leqslant 2\lVert A\rVert \epsilon.
			$$
			\item (Resonant case) if there exists $m_\ast \in \Z^d$ with $0<|m_\ast| \leqslant N$ such that
			$$
			\lvert 2 \rho-\langle m_\ast,\alpha \rangle \rvert < {\epsilon}^{\sigma},
			$$
			then
			\begin{flalign*}
				\lvert B \rvert_{r'} &\leqslant 8\left( \frac{\lVert A\rVert}{\kappa}\right) ^{\frac{1}{2}}\left( \frac{2}{r-r'} \lvert \ln \epsilon \rvert\right) ^{\frac{\tau}{2}}\times\epsilon^{\frac{-r'}{r-r'}},\\ \lVert B\rVert_0 &\leqslant 8\left( \frac{\lVert A\rVert}{\kappa}\right) ^{\frac{1}{2}}\left( \frac{2}{r-r'} \lvert \ln \epsilon \rvert\right) ^{\frac{\tau}{2}},\\
				\lvert f_{+}\rvert_{r'}&\leqslant \frac{2^{5+\tau}\lVert A\rVert\lvert \ln\epsilon\rvert^{\tau}}{\kappa(r-r')^{\tau}}\epsilon e^{-N'(r-r')}(N')^de^{Nr'}\ll \epsilon^{100}, \, N'> 2N^2.
			\end{flalign*}
			Moreover, $A_+=e^{A''}$ with $\lVert A''\rVert \leqslant 2\epsilon^{\sigma}$, $A''\in sl(2,\R)$. More accurately, we have
			$$
			MA''M^{-1}=\begin{pmatrix} it & v\\ \bar{v} & -it \end{pmatrix}
			$$
			with $\lvert t\rvert \leqslant \epsilon^{\sigma}$ and
			$$
			\lvert v \rvert\leqslant \frac{2^{4+\tau}\lVert A\rVert\lvert \ln\epsilon\rvert^{\tau}}{\kappa(r-r')^{\tau}}\epsilon e^{-\lvert m_{\ast}\rvert r}.
			$$
			
		\end{itemize}
	\end{Proposition}
	
	\begin{pf}
		We will only prove estimates for the non-resonant case because it is more delicate and the proof of the resonant case is the same compared with those in \cite{cai2022absolutely}.
		
		Let us recall that $sl(2,\R)$ is isomorphic to $su(1,1)$, which is a Lie algebra consisting of matrices of the form 
		$$
		\begin{pmatrix} 
			it & v\\ 
			\bar{v} & -it 
		\end{pmatrix}
	    $$
	    with $t \in \R, v\in \C$. The isomorphism between them is given by the map $A \to MAM^{-1}$, where
	    $$M=\frac{1}{1+i}
	    \begin{pmatrix} 
	    	1 & -i\\ 
	    	1 & i 
	    \end{pmatrix}
	    $$
	    and a straightforward calculation yields
	    $$
	    M\begin{pmatrix} 
	    	x & y+z\\ 
	    	y-z & -x 
	    \end{pmatrix}M^{-1}=
	    \begin{pmatrix} 
	    	iz & x-iy\\ 
	    	x+iy & -iz 
	    \end{pmatrix}
	    $$
	    where $x,y,z \in \R$. $SU(1,1)$ is the corresponding Lie group of $su(1,1)$. We will prove this theorem within $SU(1,1)$, which is isomorphic to $SL(2,\R)$.
	    
	    We consider the non-resonant case as follows:
	    
	    For $0<|m| \leqslant N = \frac{2}{r-r'} \lvert \ln \epsilon \rvert$, we have 
	    \begin{equation}\label{sdpc1}
	    	\lvert 2 \rho-\langle m,\alpha \rangle \rvert \geqslant {\epsilon}^{\sigma},
	    \end{equation}
	    by $(\ref{estrf})$ and $D>\frac{2}{\sigma}$, we get 
	    \begin{equation}\label{sdpc2}
	    	\left \lvert \la m,\alpha\ra \right \rvert \geqslant \frac{\kappa}{\left \lvert m \right \rvert ^{\tau}}\geqslant \frac{\kappa}{\left \lvert N \right \rvert ^{\tau}}\geqslant \epsilon^{\frac{\sigma}{2}}.
	    \end{equation}
	    
	    Indeed, it is well known that conditions $(\ref{sdpc1})$ and $(\ref{sdpc2})$ play a role in addressing the small denominator problem in KAM theory.
	    
	    We now define $g \in C^\omega _r(\T^d,su(1,1))$ such that
	    $$
	    g(\theta)=\sum_{n\in \Z^{d},0<\lvert n \rvert \leqslant N}\hat{f}(n)e^{2\pi i\la n,\theta\ra} 
	    $$
	    	    
	    From Schur's Theorem, we might as well assume $A=\begin{pmatrix}
	    	e^{i\rho} & p \\
	    	0 & e^{-i\rho} \\
	    \end{pmatrix} \in SU(1,1).$
	    The condition $||A||$ satisfying bounded gives a bound for $p$ and that is the only reason for this condition.
	    
	    Now we want to solve the cohomological equation
	    \begin{equation}\label{coequ}
	    Y(\theta + \alpha)A - AY(\theta) = -Ag(\theta).
	    \end{equation}
	    Let
	    $$
	    Y=\begin{pmatrix}
	    	y_1 & y_2 \\
	    	y_3 & y_4 \\
	    \end{pmatrix},
	    g=\begin{pmatrix}
	    	g_1 & g_2 \\
	    	g_3 & g_4 \\
	    \end{pmatrix},
	    $$
	    then we can obtain the equations
	    \begin{align*}
	    	\begin{cases}
	    		e^{i\rho}y_3(\theta + \alpha)-e^{-i\rho}y_3(\theta) = -e^{-i\rho}g_3(\theta), \\
	    		e^{i\rho}y_1(\theta + \alpha)-e^{i\rho}y_1(\theta)-py_3(\theta) = -e^{i\rho}g_1(\theta) - pg_3(\theta), \\
	    		py_3(\theta + \alpha) + e^{-i\rho}y_4(\theta + \alpha)-e^{-i\rho}y_4(\theta) = -e^{-i\rho}g_4(\theta), \\
	    		py_1(\theta + \alpha) + e^{-i\rho}y_2(\theta + \alpha)-e^{i\rho}y_2(\theta)-py_4(\theta) = -e^{i\rho}g_2(\theta)-pg_4(\theta). 
	    	\end{cases}
	    \end{align*}
	    These equations can be solved by Fourier transform. Compare the corresponding Fourier coefficients of the two sides and this shows the existence of $Y$. Apply $(\ref{sdpc1})$ twice to solve the off-diagonal and apply $(\ref{sdpc2})$ once to solve the diagonal, we can get 
	    \begin{equation}\label{estY}
	    	\lvert Y \rvert_{r'} \leqslant c\epsilon^{-3\sigma} \lvert g \rvert_r, \ 0<r'<r,
	    \end{equation}
	    where the constant only depends on $\kappa, \tau$. 	
	    
	    By the cohomological equation $(\ref{coequ})$, we obtain  
	    $$
	    Y(\theta+\alpha)=AY(\theta)A^{-1}-Ag(\theta)A^{-1}.
	    $$
	    Then we can get 
	    \begin{flalign}
	    	\notag &e^{Y(\theta+\alpha)}(Ae^{f(\theta)})e^{-Y(\theta)}\\
	    	\notag =& e^{AY(\theta)A^{-1}-Ag(\theta)A^{-1}}(Ae^{f(\theta)})e^{-Y(\theta)}\\
	    	\notag =& Ae^{\hat{f}(0)-(\mathcal{T}_Nf)(\theta)+Y(\theta)}e^{f(\theta)}e^{-Y(\theta)}\\
	    	\notag =& A \left[e^{\hat{f}(0)}+\mathcal{O}(f(\theta)-(\mathcal{T}_Nf)(\theta)+f(\theta)Y(\theta)) \right]\\
	    	\notag =& Ae^{\hat{f}(0)} \left[Id+e^{-\hat{f}(0)} \mathcal{O}((\mathcal{R}_Nf)(\theta)+f(\theta)Y(\theta)) \right]\\
	    	\notag =& A_+e^{f_+(\theta)}.
	    \end{flalign}
	    Here $\mathcal{T}_N$ is the truncation operator such that
	    $$
	    (\mathcal{T}_Nf)(\theta)=\sum_{n\in \Z^{d},\lvert n \rvert \leqslant N}\hat{f}(n)e^{2\pi i\la n,\theta\ra}
	    $$	    
	    and
	    $$
	    (\mathcal{R}_Nf)(\theta)=\sum_{n\in \Z^{d},\lvert n \rvert>N}\hat{f}(n)e^{2\pi i\la n,\theta\ra}.
	    $$	
	    
	    Therefore we define 
	    \begin{align*}
	    	\begin{cases}
	    		B(\theta)=e^{Y(\theta)}, \\
	    		A_+=Ae^{\hat{f}(0)}, \\
	    		f_+=\mathcal{O}((\mathcal{R}_Nf)(\theta)+f(\theta)Y(\theta)).
	    	\end{cases}
	    \end{align*}
    Note that although we write $D>\frac{2}{\sigma}$ i.e. $\sigma>\frac{2}{D}$, we consider $D$ to be very close to $\frac{2}{\sigma}$ in the actual process, such as $\frac{7}{D}>3\sigma>\frac{6}{D}$. Now we have estimates
    \begin{flalign*}
    &\lvert g \rvert _r \leqslant \sum_{n\in \Z^{d},0<\lvert n \rvert \leqslant N} \lvert \hat{f}(n)e^{2\pi i\la n,\theta\ra} \rvert \leqslant c\epsilon N^{d-1} < \epsilon N^{\tau-1} < \epsilon^{1-\frac{1}{D}},\\
    &\lvert Y \rvert _{r'} \leqslant c\epsilon^{-3\sigma}\lvert g \rvert _r < \epsilon^{-\frac{7}{D}} \cdot \epsilon^{1-\frac{1}{D}} < \epsilon^{1-\frac{8}{D}},\\
    &\lvert B-Id \rvert _{r'} \leqslant \lvert Y \rvert _{r'} < \epsilon^{1-\frac{8}{D}}.
    \end{flalign*}
    \begin{flalign*}
    \lvert (\mathcal{R}_Nf)(\theta)\rvert_{r'} & \leqslant \sum_{n\in \Z^{d},\lvert n \rvert>N}\lvert \hat{f}(n)e^{2\pi i\la n,\theta\ra}\rvert_{r'}\\
    & \leqslant c\lvert f \rvert _r e^{-N(r-r')}(N+\frac{1}{r-r'})^{d}\\
    & \leqslant c\epsilon e^{-2\log \frac{1}{\epsilon}}(N)^{d}\\
    & < \epsilon \cdot \epsilon^{2}\cdot \epsilon^{-\frac{1}{D}}\\
    & <\epsilon^{3-\frac{1}{D}}.
    \end{flalign*}
    \begin{flalign*}
    & \lvert f_+ \rvert _{r'} \leqslant \lvert fY \rvert _{r'} + \lvert (\mathcal{R}_Nf)\rvert_{r'} \leqslant c\epsilon \cdot \epsilon^{1-\frac{8}{D}} + \epsilon^{3-\frac{1}{D}} < \epsilon^{2-\frac{8}{D}}.\\
    & \lVert A_+-A\rVert\leqslant \lVert A\rVert \lVert Id-e^{\hat{f}(0)} \rVert \leqslant 2\lVert A\rVert \epsilon.
    \end{flalign*}
    This finish the proof of Propositon $\ref{pro3.1}$.
	\end{pf}
	
	\begin{Remark}
		This version of analytic KAM theorem is a perfect combination of Eliasson \cite{MR1167299} and Hou-You \cite{hou2012almost}. While the resonant case which absorbs the essence of Hou-You stays the same, the refined non-resonant case avoids eliminating irrelevant non-resonant terms via Eliasson's way compared with Cai \cite{cai2022absolutely}. This essentially reduces the norm of conjugation maps, which is later vital for us to ensure that the final loss of regularity is  independent of the initial $k$ for our $C^k$ reducibility theorem. Note that this is at all not considered in Cai \cite{cai2022absolutely} as almost reducibility does not really care about the convergence of the conjugation maps, but reducibility does. 
	\end{Remark}

	\subsection{$C^k$ reducibility}
	As planned, we will present the quantitative $C^k$ reducibility via analytic approximation \cite{zehnder1975generalized}.
	
	The crucial improvement here compared with \cite{cai2022absolutely} lies at the point where we are able to obtain reducibility results in a KAM step. In fact, when the resonant steps are well-separated from each other, it allows for improved control over the conjugation maps during the inductive argument. And then, after we apply certain conditions to the rotation number of cocycle $(\alpha, Ae^{f(\theta)})$, we can find that there are only finite resonant steps in the whole KAM iteration process. This, in return, will give us the best possible initial regularity $k$ through the technique.
	
	Before we prove the main theorem, we will first cite the Proposition 3.1 of \cite{cai2022absolutely} to simplify the main proof. We first recall some notations given in \cite{cai2022absolutely}.
	
	Let $\{f_{j}\}_{j\geqslant 1}$, $f_{j}\in C_{\frac{1}{j}}^{\omega}(\T^{d},sl(2,\R))$ be the analytic sequence approximating $f\in C^k(\T^{d},sl(2,\R))$ which satisfies $(\ref{aa})$.
	
	For $0<r'<r$, denote
	\begin{equation}\label{denote1}
		\epsilon_0'(r,r')=\frac{c}{(2\lVert A\rVert)^{\tilde{D}}}(r-r')^{D\tau},
	\end{equation}
	and
	\begin{equation}\label{denote2}
		\epsilon_m=\frac{c}{(2\lVert A\rVert)^{\tilde{D}}m^{D\tau+\frac{1}{2}}}, m\in\Z^+.
	\end{equation}
    where $c$ depends on $\kappa, \tau, d$, $D,\tilde{D} \in \Z$ depend on $\sigma$.
	
	Then for any $0<s\leqslant\frac{1}{6D\tau+3}$ fixed, there exists $m_0$ such that for any $m\geqslant m_0$ we can get 
	\begin{equation}\label{denote3}
		\frac{c}{(2\lVert A\rVert)^{\tilde{D}}m^{D\tau+\frac{1}{2}}}\leqslant \epsilon_0'\left( \frac{1}{m},\frac{1}{m^{1+s}}\right) .
	\end{equation}
	
	We will begin with $M> \max\{\frac{(2\lVert A\rVert)^{\tilde{D}}}{c},m_0\}$, $M\in\N^+$. Let $l_j=M^{(1+s)^{j-1}}$, $j\in \N^+$. Because $l_j$ is not an integer, we pick $[l_j]+1$ instead of $l_j$.
	
	Now, let $\Omega=\{n_1,n_2,n_3,\cdots\}$ denote the sequence of all resonant steps. In other words, the $(n_j)$-th step is obtained by resonant case. By analytic approximation $(\ref{aa})$ and Proposition $\ref{pro3.1}$ in each iteration step, we can establish the following almost reducibility result concerning each $(\alpha, Ae^{f_{l_j}(\theta)})$ by applying induction.
	
	\begin{Proposition}\cite{cai2022absolutely} \cite{cai2022reducibility}\label{pro3.2}
		Let $\alpha\in {\rm DC}_d(\kappa,\tau)$, $\sigma<\frac{1}{6}$. Assume that $A\in SL(2,\R)$ satisfying $A$ abounded, $f\in C^k(\T^d, sl(2,\R))$ with $k>(D+2)\tau+2$ and $\{f_j\}_{j\geqslant 1}$ be as Sect. 2.4. There exists $\bar \epsilon=\bar \epsilon(\kappa,\tau,d,k, \lVert A \rVert,\sigma)$ such that if 
		\begin{flalign}\label{estkf}
		\lVert f\rVert_k \leqslant \bar \epsilon \leqslant \epsilon_0'\left( \frac{1}{l_1},\frac{1}{l_2}\right) ,
	    \end{flalign} 
		then there exist $B_{l_j}\in C^{\omega}_{\frac{1}{l_{j+1}}}(2\T^d, SL(2,\R))$, $A_{l_j}\in SL(2,\R)$ and $f^{'}_{l_j}\in C^{\omega}_{\frac{1}{l_{j+1}}}(\T^d, sl(2,\R))$ such that
		\begin{flalign}\label{conj1}
		B_{l_j}(\theta+\alpha)(Ae^{f_{l_j}(\theta)})B^{-1}_{l_j}(\theta)=A_{l_j}e^{f_{l_j}^{'}(\theta)},
	    \end{flalign}
		with estimates
		\begin{flalign}
			\label{estrb1} & \lvert B_{l_j}(\theta)\rvert_{\frac{1}{l_{j+1}}} \leqslant 64\left(\frac{\lVert A\rVert}{\kappa}\right) \left( \frac{2}{\frac{1}{l_j}-\frac{1}{l_{j+1}}}  \ln {\frac{1}{\epsilon_{l_j}}}\right) ^{\tau}\times{\epsilon_{l_j}}^{-\frac{\frac{2}{l_{j+1}}}{\frac{1}{l_j}-\frac{1}{l_{j+1}}}} \leqslant \epsilon_{l_{j}}^{-\frac{\sigma}{2}-s},\\
			\label{est0b1} & \lVert B_{l_j}(\theta)\rVert_0 \leqslant 64\left( \frac{\lVert A\rVert}{\kappa}\right) \left( \frac{2}{\frac{1}{l_j}-\frac{1}{l_{j+1}}}  \ln {\frac{1}{\epsilon_{l_j}}}\right) ^{\tau} \leqslant \epsilon_{l_{j}}^{-\frac{\sigma}{2}},\\
			\label{estdeg1} & \lvert \deg B_{l_j}\rvert \leqslant 4l_j \ln {\frac{1}{\epsilon_{l_j}}},\\
			\label{estrf'} & \lvert f_{l_j}^{'}(\theta)\rvert_{\frac{1}{l_{j+1}}} \leqslant \epsilon_{l_{j}}^{2-\frac{8}{D}},\ \ \lVert A_{l_j}\rVert\leqslant 2\lVert A\rVert.
		\end{flalign}
		Moreover, there exists $\overline f_{l_j}\in C^{k_0}(\T^d, sl(2,\R))$ with $k_0 \in \N, k_0 \leqslant \frac{k-10\tau-3}{1+s}$ such that
		\begin{flalign}\label{conj2}
		B_{l_j}(\theta+\alpha)(Ae^{f(\theta)})B^{-1}_{l_j}(\theta)=A_{l_j}e^{\overline f_{l_j}(\theta)},
	    \end{flalign}
		with estimate
		\begin{flalign}\label{estkf-}
		\lVert \overline f_{l_j}(\theta)\rVert_{k_0}\leqslant \epsilon_{l_{j}}^{\frac{3}{D}}.
	    \end{flalign}		
	\end{Proposition}
	
	\begin{Remark}\label{rem3.2}
		By the new analytic KAM scheme and the original proof method in  \cite{cai2022absolutely} \cite{cai2022reducibility}, we can also get the results of the proposition. The fact that the existing estimates $(\ref{estrf'})$ and $(\ref{estkf-})$ are more delicate than the original is also due to the use of the new $C^\omega$ KAM theorem. 
	\end{Remark}
	
	\begin{Remark}\label{rem3.3}
		The estimates of $(\ref{estdeg1})$ can be done better in the non-resonant case because the non-resonant step does not change the degree, thus the estimates of the degree in the non-resonant case can be considered as the estimates of the resonant step closest to this step.
	\end{Remark}
		
	With Proposition \ref{pro3.2} in hand, we are going to transfer all the estimates from $(\alpha, Ae^{f_{l_j}(\theta)})$ to $(\alpha, Ae^{f(\theta)})$ through analytic approximation. We establish the following quantitative $C^k$ reducibility theorem.
	
	\begin{Theorem}\label{thm3.1}
		Let $\alpha \in {\rm DC}_d(\kappa,\tau)$, $\sigma<\frac{1}{6}$, $A\in SL(2,\R)$ and $f\in C^k(\T^{d}, \\ sl(2,\R))$ with $k>(D+2)\tau+2$. Then there exists $\epsilon_0=\epsilon_0(\kappa,\tau,d,k,\lVert A\rVert,\sigma)$ such that if 
		\begin{flalign}\label{estkf2}
			\lVert f\rVert_k \leqslant \epsilon_0 \leqslant \epsilon_0'\left( \frac{1}{l_1},\frac{1}{l_2}\right) 
		\end{flalign}
		and 
		\begin{itemize}
			\item if $\rho(\alpha, Ae^f)$ is Diophantine with respect to $\alpha $: $\rho(\alpha,Ae^f) \in {\rm DC}_d ^\alpha(\gamma,\tau)$, then there exists two constants $C_1=C_1(\gamma,\kappa,\tau,d,k,||A||,\sigma)$, $C_2=C_2(\gamma,\kappa,\tau,d,k,||A||,\sigma)$ and $B_1\in C^{k_0}(2\T^d, SL(2,\R))$ with $k_0 \in \N$, $k_0\leqslant \frac{k-10\tau-3}{1+s}$ such that
			\begin{flalign}\label{trans1}
				B_1(\theta+\alpha)(Ae^{f(\theta)})B^{-1}_1(\theta)=R_\phi \in SL(2,\R), \ \phi \notin \Z, 
			\end{flalign}
			with estimates
			\begin{flalign}\label{estb}
				\lVert B_1 \rVert_{k_0} \leqslant C_1, \ |\deg B_1 | \leqslant C_2. 
			\end{flalign}
			\item if $\rho(\alpha,Ae^f)$ is rational with respect to $\alpha $: $2\rho(\alpha,Ae^f) =\la m_0,\alpha  \ra \mod \Z$ for some $m_0\in \Z^d$, then there exists $B_2\in C^{k_0}(2\T^d, SL(2,\R))$ with $k_0 \in \N$, $k_0\leqslant \frac{k-10\tau-3}{1+s}$ such that
			\begin{flalign}\label{trans2}
				B_2(\theta+\alpha)(Ae^{f(\theta)})B^{-1}_2(\theta)=\tilde A_2 \in SL(2,\R), 
			\end{flalign}
		    with
		    \begin{flalign}\label{estrot}
			\rho(\alpha, \tilde A_2)=0.
		    \end{flalign}
		\end{itemize}
	\end{Theorem}
	
	\begin{Remark}\label{rem3.4}
		In the inequality $k_0\leqslant \frac{k-10\tau-3}{1+s}$, due to the precise choice of $k, \tau, s$, we can also pick $k_0 = [k-10\tau-3]$. Here ``$[x]$'' stands for the integer part of $x$.
	\end{Remark}
	
	\begin{pf}
		\textbf{(Diophantine case)}
		By $(\ref{estkf2})$ and Proposition $\ref{pro3.2}$, take $\epsilon_0 = \bar \epsilon$ as in Proposition $\ref{pro3.2}$ and apply it to cocycle $(\alpha,Ae^{f(\theta)})$, then there exists $B_{l_j}\in C^{\omega}_{\frac{1}{l_{j+1}}}(2\T^d, SL(2,\R))$, $A_{l_j}\in SL(2,\R)$ and $\overline f_{l_j}\in C^{k_0}(\T^d, sl(2,\R))$ with $k_0 \in \N$, $k_0 \leqslant \frac{k-10\tau-3}{1+s}$ such that
		\begin{flalign}\label{conj6}
		     B_{l_j}(\theta+\alpha)(Ae^{f(\theta)})B^{-1}_{l_j}(\theta)=A_{l_j}e^{\overline f_{l_j}(\theta)},
	    \end{flalign}
		with estimates
		     \begin{flalign}\label{estb1}
		     	& \lvert B_{l_j}(\theta)\rvert_{\frac{1}{l_{j+1}}}\leqslant \epsilon_{l_{j}}^{-\frac{\sigma}{2}-s}, \ \lVert B_{l_j}(\theta)\rVert_0\leqslant \epsilon_{l_{j}}^{-\frac{\sigma}{2}}, \ |\deg B_{l_j}| \leqslant 4l_j \ln \frac{1}{\epsilon_{l_{j}}}, \\
			    \label{estkf-2} & \lVert A_{l_j} \rVert \leqslant 2\lVert A \rVert,\ \ \lVert \overline f_{l_j}(\theta) \rVert_{k_0}\leqslant \epsilon_{l_{j}}^{\frac{3}{D}}.
		    \end{flalign}
	
		In the last part of analytic approximation $(\ref{aa})$, taking a telescoping sum from $j$ to $+\infty$, we get
		\begin{equation}\label{aa1}
			\lVert f(\theta)-f_{l_j}(\theta)\rVert_0\leqslant \frac{c}{(2\lVert A\rVert)^{\tilde{D}}l_1^{D\tau+\frac{1}{2}}l_{j}^{k-1}},
		\end{equation}
		\begin{equation}\label{aa2}
			\lVert f(\theta)\rVert_0+\lVert f_{l_j}(\theta)\rVert_0\leqslant \frac{c}{(2\lVert A\rVert)^{\tilde{D}}M^{D\tau+\frac{1}{2}}}+\frac{\tilde c}{(2\lVert A\rVert)^{\tilde{D}}M^{D\tau+\frac{1}{2}}}.
		\end{equation}
		Thus by $(\ref{conj6})-(\ref{aa2})$, we have
		\begin{flalign}\label{est0f-}
			\begin{split}
			\lVert \overline f_{l_j}(\theta)\rVert_0 \leqslant & \lVert f_{l_j}^{'}(\theta)\rVert_0 + \lVert A^{-1}_{l_j}\rVert \lVert B_{l_j}(\theta+\alpha)(Ae^{f(\theta)}-Ae^{f_{l_j}(\theta)})B^{-1}_{l_j}(\theta)\rVert_0 \\
			\leqslant & \epsilon_{l_j}^{2-\frac{8}{D}}+2\lVert A\rVert \times \epsilon_{l_j}^{-\sigma} \times \frac{c}{(2\lVert A\rVert)^{\tilde{D}}l_1^{D \tau + \frac{1}{2}}l_j^{k-1}} \\ 
			\leqslant & \epsilon_{l_j}^{1+s}.
			\end{split}
		\end{flalign}
		
		Since $\rho(\alpha,Ae^f) \in {\rm DC}_d ^\alpha(\gamma,\tau)$, for any $m \in \Z^d$, we have
		\begin{flalign*}
			&\lVert 2\rho(\alpha,A_{l_j}e^{\overline f_{l_j}(\theta)})-\la m,\alpha \ra \rVert_{\R / \Z} \\
			= &\lVert 2\rho(\alpha,Ae^{f(\theta)}) - \la \deg B_{l_j}, \alpha \ra -\la m,\alpha \ra \rVert_{\R / \Z} \\
			\geqslant &\frac{\gamma}{(|m+\deg B_{l_j}|+1)^\tau} \\
			\geqslant &\frac{\gamma (1+|\deg B_{l_j}|)^{-\tau}}{(|m|+1)^\tau},
		\end{flalign*}
		which implies $\rho(\alpha,A_{l_j}e^{\overline f_{l_j}(\theta)}) \in {\rm DC}_d ^\alpha(\gamma (1+|\deg B_{l_j}|)^{-\tau},\tau)$.
		
		By $(\ref{estb1})$, we can obtain 
		$$
		\epsilon_{l_j}^\frac{(1+s) \sigma}{2} (1+|\deg B_{l_j}|)^\tau \leqslant \left[ \frac{c}{(2\lVert A\rVert)^{\tilde{D}}{l_j}^{D\tau+\frac{1}{2}}} \right]^\frac{(1+s) \sigma}{2} \left( 1+4l_j \ln \frac{1}{\epsilon_{l_j}}\right) ^\tau.
		$$
		Obviously, let $j \to \infty$ and then the right-hand side of the inequality to zero. Therefore, for any given $\gamma > 0$, there exists $j'=j'(\gamma) \in \N$ such that
		$$
		\epsilon_{l_{j'}}^\frac{(1+s) \sigma}{2} (1+|\deg B_{l_{j'}}|)^\tau \leqslant \gamma
		$$
		and then
        \begin{equation}\label{estfz}
		   \gamma (1+|\deg B_{l_{j'}}|)^{-\tau} \geqslant \epsilon_{l_{j'}}^\frac{(1+s) \sigma}{2} = \epsilon_{l_{j'+1}}^\frac{\sigma}{2}.
	    \end{equation}
        
        On the other hand, for $0<|m| \leqslant N_{l_{j'+1}}=\frac{2}{\frac{1}{l_{j'+1}}-\frac{1}{l_{j'+2}}} \lvert \ln \epsilon_{l_{j'+1}} \rvert$, we have
        \begin{equation}\label{estfm}
        	\frac{1}{(|m|+1)^\tau} \geqslant \frac{1}{\left(\frac{2}{\frac{1}{l_{j'+1}}-\frac{1}{l_{j'+2}}} \ln \frac{1}{\epsilon_{l_{j'+1}}}+1\right)^ \tau} \geqslant 2\epsilon_{l_{j'+1}}^\frac{\sigma}{2}.
        \end{equation}
        Then by $(\ref{estfz})$, $(\ref{estfm})$, we can get 
        \begin{flalign*}
        	&\lVert 2\rho(\alpha,A_{l_{j'}})-\la m,\alpha \ra \rVert_{\R / \Z} \\
        	\geqslant &\lVert 2\rho(\alpha,A_{l_{j'}}e^{\overline f_{l_{j'}}(\theta)})-\la m,\alpha \ra \rVert_{\R / \Z}-| 2\rho(\alpha,A_{l_{j'}}e^{\overline f_{l_{j'}}(\theta)})-2\rho(\alpha,A_{l_{j'}})| \\
        	\geqslant &\frac{\gamma (1+|\deg B_{l_{j'}}|)^{-\tau}}{(|m|+1)^\tau}-2\epsilon_{l_{j'}}^\frac{1+s}{2} \\
        	\geqslant &\epsilon_{l_{j'+1}}^\frac{\sigma}{2} \cdot 2\epsilon_{l_{j'+1}}^\frac{\sigma}{2} - 2\epsilon_{l_{j'+1}}^\frac{1}{2} \\
        	\geqslant &\epsilon_{l_{j'+1}}^\sigma,
        \end{flalign*}
        which means the $(j'+1)$-th step is non-resonant with
        $$
        B_{l_{j'+1}}=\widetilde{B}_{l_{j'}} \circ B_{l_{j'}}, \ \lvert \widetilde B_{l_{j'}}(\theta)- Id\rvert_{\frac{1}{l_{j'+2}}} \leqslant \epsilon_{l_{j'+1}}^{1-\frac{8}{D}}, \ \deg B_{l_{j'+1}}=\deg B_{l_{j'}}.
        $$
        
        Assume that for $l_n$, $j'+1 \leqslant n \leqslant j_0$, we have 
        $$
        B_{l_n}(\theta+\alpha)(Ae^{f_{l_{n}}(\theta)})B^{-1}_{l_n}(\theta)=A_{l_n}e^{f_{l_n}'(\theta)},
        $$
        which is equivalent to
        $$
        B_{l_n}(\theta+\alpha)(Ae^{f(\theta)})B^{-1}_{l_n}(\theta)=A_{l_n}e^{f_{l_n}'(\theta)}+B_{l_n}(\theta+\alpha)(Ae^{f(\theta)}-Ae^{f_{l_{n}(\theta)}})B^{-1}_{l_n}(\theta),
        $$
        rewrite that
        \begin{equation}\label{conj7}
            B_{l_n}(\theta+\alpha)(Ae^{f(\theta)})B^{-1}_{l_n}(\theta)=A_{l_n}e^{\overline f_{l_n}(\theta)},
        \end{equation}
        with estimates
        \begin{equation}\label{estb2}
        	B_{l_{n}}=\widetilde{B}_{l_{n-1}} \circ B_{l_{n-1}}, \ \lvert \widetilde B_{l_{n-1}}(\theta)- Id\rvert_{\frac{1}{l_{n+1}}} \leqslant \epsilon_{l_n}^{1-\frac{8}{D}}, \ \deg B_{l_{n-1}}=\deg B_{l_{j'}}.
        \end{equation}
        
        Therefore, by $(\ref{estb2})$ and the Diophantine condition on $\rho(\alpha,Ae^f)$, for $0<|m| \leqslant \frac{2}{\frac{1}{l_{n+1}}-\frac{1}{l_{n+2}}} \lvert \ln \epsilon_{l_{n+1}} \rvert$, we obtain
        \begin{flalign*}
         	&\lVert 2\rho(\alpha,A_{l_n})-\la m,\alpha \ra \rVert_{\R / \Z} \\
         	\geqslant &\lVert 2\rho(\alpha,A_{l_n}e^{\overline f_{l_n}(\theta)})-\la m,\alpha \ra \rVert_{\R / \Z}-| 2\rho(\alpha,A_{l_n}e^{\overline f_{l_n}(\theta)})-2\rho(\alpha,A_{l_n})| \\
         	\geqslant &\frac{\gamma (1+|\deg B_{l_{j'}}|)^{-\tau}}{(|m|+1)^\tau}-2\epsilon_{l_n}^\frac{1+s}{2} \\
         	\geqslant &\epsilon_{l_{j'+1}}^\frac{\sigma}{2} \cdot 2\epsilon_{l_{n+1}}^\frac{\sigma}{2} - 2\epsilon_{l_{n+1}}^\frac{1}{2} \\
         	\geqslant &\epsilon_{l_{n+1}}^\sigma, \ \forall j'+1 \leqslant n \leqslant j_0.
        \end{flalign*}
        This means the $(j_0 + 1)$-th step is still non-resonant with estimates
        $$
        B_{l_{j_0 + 1}}=\widetilde{B}_{l_{j_0}} \circ B_{l_{j_0}}, \ \lvert \widetilde B_{l_{j_0}}(\theta)- Id\rvert_{\frac{1}{l_{j_0 + 2}}} \leqslant \epsilon_{l_{j_0 + 1}}^{1-\frac{8}{D}}, \ \deg B_{l_{j_0}}=\deg B_{l_{j'}}.
        $$
        
        In conclusion, we know that there are at most finitely many resonant steps in the iteration process. Assume that $n_q$ is the last resonant step, then for $\forall j > n_q$, we have 
        $$
        B_{l_j}(\theta+\alpha)(Ae^{f_{l_{j}}(\theta)})B^{-1}_{l_j}(\theta)=A_{l_j}e^{f_{l_j}^{'}(\theta)},
        $$
        with estimates
        \begin{equation}\label{estb3}
        	B_{l_{j}}=\widetilde{B}_{l_{j-1}} \circ B_{l_{j-1}}, \ \lvert \widetilde B_{l_{j-1}}(\theta)- Id\rvert_{\frac{1}{l_{j+1}}} \leqslant \epsilon_{l_j}^{1-\frac{8}{D}}, \ \deg B_{l_{j-1}}=\deg B_{l_{n_q}}.
        \end{equation}
        
        Denote $B_1=\lim \limits_{j\to\infty}B_{l_j}$, $\tilde A_1=\lim \limits_{j\to\infty}A_{l_j} \in SL(2,\R)$. Notice that $\rho(\alpha, \tilde A_1) \neq 0$, otherwise it will contradict to $\rho(\alpha,Ae^f) \in {\rm DC}_d ^\alpha(\gamma,\tau)$. Thus $\tilde A_1$ can only be standard rotation in $SL(2,\R)$, which is the case of $(\ref{trans1})$.
        
        By $(\ref{estb3})$ and Cauchy estimates, for $\forall j > n_q$, we can calculate
        \begin{flalign*}
        	\lVert \widetilde B_{l_{j}} - Id \rVert_{k_0} \leqslant & \sup_{\substack{\lvert l\rvert\leqslant k_0, \theta \in \T^{d}}}\lVert (\partial^{l_1}_{\theta_1} \cdots \partial^{l_d}_{\theta_d}) (\widetilde B_{l_{j}} - Id) \rVert \\ 
        	\leqslant & (k_0)!(l_{j+2})^{k_0}\lvert \widetilde B_{l_{j}} - Id \rvert_{\frac{1}{l_{j+2}}} \\
        	\leqslant & (k_0)!(l_{j+1})^{(1+s)k_0} \epsilon_{l_{j+1}}^{1-\frac{8}{D}} \\
        	\leqslant & \frac{C}{l_{j+1}^{(1-\frac{8}{D})(D\tau+\frac{1}{2})-(1+s)k_0}},
        \end{flalign*} 
        where $C$ does not depend on $j$. 
        
        Note that while we write the assumption as $k>(D+2)\tau+2$, in the actual operational process we simply choose $k=[(D+2)\tau+2]+1$. Therefore, the value of $k$ is entirely determined by the parameter $D$.  
        
        Thus if we pick $k_0 \leqslant \frac{k-10\tau-3}{1+s}$, we have 
        $$
        \lVert \widetilde B_{l_{j}} \rVert_{k_0} \leqslant 1+ \frac{C}{l_{j+1}^{\frac{1}{6}}}.
        $$  
        Since we pick $l_1=M$ sufficiently large, then
        \begin{flalign*}
        	\lVert B_1 \rVert_{k_0} \leqslant & \lVert \prod\limits_{j=n_q+1}^{\infty} \widetilde B_{l_{j}} \rVert_{k_0} \lVert B_{l_{n_q}} \rVert_{k_0} \\
        	\leqslant & \left( \prod\limits_{j=n_q+1}^{\infty} \left( 1+ \frac{C}{l_{j+1}^{\frac{1}{6}}}\right) \right) \sup_{\substack{\lvert l\rvert\leqslant k_0, \theta \in \T^{d}}}\lVert (\partial^{l_1}_{\theta_1} \cdots \partial^{l_d}_{\theta_d}) (B_{l_{n_q}}(\theta)) \rVert \\ 
        	\leqslant & 2(k_0)!(l_{n_q+1})^{k_0}\lvert B_{l_{n_q}}\rvert_{\frac{1}{l_{n_q+1}}} \\
        	\leqslant & {l_{n_q}}^{(\frac{\sigma}{2}+s)(D\tau+\frac{1}{2})+(1+s)k_0}.
        \end{flalign*}            
        The estimate of $\deg B_1$ is clearly valid, which gives $(\ref{estb})$.
        
        \textbf{(Rational case)}
        By $(\ref{estkf2})$ and Proposition $\ref{pro3.2}$, take $\epsilon_0 = \bar \epsilon$ as in Proposition $\ref{pro3.2}$ and apply it to cocycle $(\alpha,Ae^{f(\theta)})$, then there exists $B_{l_j}\in C^{\omega}_{\frac{1}{l_{j+1}}}(2\T^d, \\  SL(2,\R))$, $A_{l_j}\in SL(2,\R)$ and $\overline f_{l_j}\in C^{k_0}(\T^d, sl(2,\R))$ with $k_0 \in \N$, $k_0\leqslant \frac{k-10\tau-3}{1+s}$ such that
        \begin{flalign}\label{conj8}
        	B_{l_j}(\theta+\alpha)(Ae^{f(\theta)})B^{-1}_{l_j}(\theta)=A_{l_j}e^{\overline f_{l_j}(\theta)},
        \end{flalign}
        with estimates
        \begin{flalign}\label{estb4}
        	& \lvert B_{l_j}(\theta)\rvert_{\frac{1}{l_{j+1}}}\leqslant \epsilon_{l_{j}}^{-\frac{\sigma}{2}-s}, \ \lVert B_{l_j}(\theta)\rVert_0\leqslant \epsilon_{l_{j}}^{-\frac{\sigma}{2}}, \ |\deg B_{l_j}| \leqslant 4l_j \ln \frac{1}{\epsilon_{l_{j}}}, \\
            \label{estkf-3} & \lVert A_{l_j} \rVert \leqslant 2\lVert A \rVert,\ \ \lVert \overline f_{l_j}(\theta) \rVert_{k_0}\leqslant \epsilon_{l_{j}}^{\frac{3}{D}}.
        \end{flalign}
        
        Since $\rho(\alpha,Ae^f) = \frac{\la m_0,\alpha \ra}{2} \mod \frac{\Z}{2}$, $m_0\in \Z^d$, we have
        \begin{flalign*}
        	\rho(\alpha,A_{l_j}e^{\overline f_{l_j}(\theta)}) & = \rho(\alpha,Ae^{f(\theta)}) - \frac{\la \deg B_{l_j},\alpha \ra}{2} \mod \frac{\Z}{2} \\
        	& = \frac{\la m_0 - \deg B_{l_j},\alpha \ra}{2} \mod \frac{\Z}{2}.
        \end{flalign*}
        From now on, we omit ``$\mod \frac{\Z}{2}$" for simplicity. By $(\ref{rot2})$ and the proof of Proposition $\ref{pro3.2}$, we know the rotation number is invariant in non-resonant case. 
        
        Let $m_{l_{n_i}} \in \Z^d, i=1,2,3, \cdots $ represent resonant sites of the $(n_i)$-th step with $0<|m_{l_{n_i}}| \leqslant N_{l_{n_i}}=\frac{2}{\frac{1}{l_{n_i}}-\frac{1}{l_{n_i+1}}} \lvert \ln \epsilon_{l_{n_i}}\rvert$. By Claim 1 in \cite{cai2022absolutely}, we have $N_{l_{n_{i+1}}} \gg 4N_{l_{n_i}}, i=1,2,3, \cdots$. So there must exist $j \in \Z$ sufficiently large, provided that there are $q-1$ resonant steps before the $j$-th step, such that
        $$
        m_0 - (m_{l_{n_1}}+m_{l_{n_2}}+ \cdots +m_{l_{n_{q-1}}}) = m',
        $$
        where $m' \in \Z^d$ with $0<|m'| \leqslant N_{l_j}=\frac{2}{\frac{1}{l_{j}}-\frac{1}{l_{j+1}}} \lvert \ln \epsilon_{l_{j}}\rvert$ and $N_{l_j} \gg 2N_{l_{n_{q-1}}} \gg  m_{l_{n_1}}+m_{l_{n_2}}+ \cdots +m_{l_{n_{q-1}}}$. And then
        \begin{equation}\label{rot0} 
        	\begin{split}
        	\rho(\alpha,A_{l_j}e^{\overline f_{l_j}(\theta)}) &= \frac{\la m_0 - \deg B_{l_j},\alpha \ra}{2} \\ 
        	&= \frac{\la m_0 - (m_{l_{n_1}}+m_{l_{n_2}}+ \cdots +m_{l_{n_{q-1}}}+m'),\alpha \ra}{2} \\
        	&= 0.
        	\end{split}
        \end{equation}
        Since all the steps between $(n_{q-1})$-th step and $j$-th step are non-resonant and by $(\ref{rot2})$, $(\ref{est0f-})$ and $(\ref{rot0})$, we have
        \begin{flalign*}
        	& \lVert 2\rho(\alpha,A_{l_{j-1}})-\la m',\alpha \ra \rVert_{\R / \Z} \\
        	\leqslant & \lVert 2\rho(\alpha,A_{l_{j-1}}e^{\overline f_{l_{j-1}}(\theta)})-\la m',\alpha \ra \rVert_{\R / \Z} + |2\rho(\alpha,A_{l_{j-1}}e^{\overline f_{l_{j-1}}(\theta)})-2\rho(\alpha,A_{l_{j-1}})| \\
        	\leqslant & |2\rho(\alpha,A_{l_{j}}e^{\overline f_{l_j}(\theta)})| + 2\epsilon_{l_{j-1}}^\frac{1+s}{2} \\
        	= & 0 + 2\epsilon_{l_{j}}^\frac{1}{2} \\
        	< & \epsilon_{l_{j}}^\sigma.
        \end{flalign*}
        Thus the $j$-th step is the $(n_q)$-th resonant step and $m'$ is the unique resonant site $m_{l_{n_q}}$ with $0<|m_{l_{n_q}}| \leqslant N_{l_{n_q}}=\frac{2}{\frac{1}{l_{n_q}}-\frac{1}{l_{n_q+1}}} \lvert \ln \epsilon_{l_{n_q}}\rvert$. 
        
        Now we apply Proposition $\ref{pro3.2}$ to cocycle $(\alpha,A_{l_{n_q}}e^{\widetilde f_{l_{n_q}}(\theta)})$, then we can get
        $$
        \widetilde{B}_{l_{n_q}}(\theta+\alpha)(A_{l_{n_q}}e^{\widetilde f_{l_{n_q}}(\theta)})\widetilde{B}_{l_{n_q}}^{-1}(\theta)=A_{l_{{n_q}+1}}e^{f_{l_{{n_q}+1}}^{'}(\theta)},
        $$
        which gives
        $$
        {B}_{l_{{n_q}+1}}(\theta+\alpha)(Ae^{{f_{l_{{n_q}+1}}}(\theta)}){B}_{l_{{n_q}+1}}^{-1}(\theta)=A_{l_{{n_q}+1}}e^{f_{l_{{n_q}+1}}^{'}(\theta)},
        $$
        where $B_{l_{n_q+1}}=\widetilde{B}_{l_{n_q}} \circ B_{l_{n_q}} \in C^{\omega}_{\frac{1}{l_{n_q+2}}}(2\T^{d},SL(2,\R))$.
        
        Since $\rho(\alpha,A_{l_{n_q}}e^{\overline f_{l_{n_q}}(\theta)}) = \rho(\alpha,Ae^{\overline f_{l_{n_q}}(\theta)}) = 0$, for $m \in \Z^d$ with $0<|m| \leqslant N_{l_{n_q+1}}=\frac{2}{\frac{1}{l_{n_q+1}}-\frac{1}{l_{n_q+2}}} \lvert \ln \epsilon_{l_{n_q+1}}\rvert$, we have 
        \begin{flalign*}
        	& \lVert 2\rho(\alpha,A_{l_{n_q}})-\la m,\alpha \ra \rVert_{\R / \Z} \\
        	\geqslant & \lVert 2\rho(\alpha,A_{l_{n_q}}e^{\overline f_{l_{n_q}}(\theta)})-\la m,\alpha \ra \rVert_{\R / \Z} - |2\rho(\alpha,A_{l_{n_q}}e^{\overline f_{l_{n_q}}(\theta)})-2\rho(\alpha,A_{l_{n_q}})| \\
        	\geqslant & \frac{\kappa}{|m|^ \tau} - 2\epsilon_{l_{n_q}}^\frac{1+s}{2} \\
        	\geqslant & \frac{\kappa}{\left( \frac{2}{\frac{1}{l_{n_q+1}}-\frac{1}{l_{n_q+2}}} \ln \frac{1}{\epsilon_{l_{n_q+1}}}\right) ^\tau} - 2\epsilon_{l_{n_q+1}}^\frac{1}{2} \\
        	\geqslant & 2\epsilon_{l_{n_q+1}}^\frac{\sigma}{2} - 2\epsilon_{l_{n_q+1}}^\frac{1}{2} \\
        	\geqslant & \epsilon_{l_{n_q+1}}^\sigma,
        \end{flalign*}
        which means the $(n_q+1)$-th step is non-resonant with
        $$
        \deg B_{l_{n_q+1}} = \deg \widetilde B_{l_{n_q}} + \deg B_{l_{n_q}} = m_0.
        $$
        
        Assume that for $l_j$, $n_q+1 \leqslant j \leqslant j_0$, we already have
        $$
        B_{l_j}(\theta+\alpha)(Ae^{f_{l_{j}}(\theta)})B^{-1}_{l_j}(\theta)=A_{l_j}e^{f_{l_j}^{'}(\theta)},
        $$
        which is equivalent to
        $$
        B_{l_j}(\theta+\alpha)(Ae^{f(\theta)})B^{-1}_{l_j}(\theta)=A_{l_j}e^{f_{l_j}^{'}(\theta)}+B_{l_j}(\theta+\alpha)(Ae^{f(\theta)}-Ae^{f_{l_{j}}})B^{-1}_{l_j}(\theta),
        $$
        rewrite that
        \begin{equation}\label{conj9}
        	B_{l_j}(\theta+\alpha)(Ae^{f(\theta)})B^{-1}_{l_j}(\theta)=A_{l_j}e^{\overline f_{l_j}(\theta)},
        \end{equation}
        with estimates
        \begin{equation}\label{estb5}
        	B_{l_{j}}=\widetilde{B}_{l_{j-1}} \circ B_{l_{j-1}}, \lvert \widetilde B_{l_{j-1}}(\theta)- Id\rvert_{\frac{1}{l_{j+1}}} \leqslant \epsilon_{l_j}^{1-\frac{8}{D}}, \deg B_{l_{j-1}}= m_0.
        \end{equation}
        
        Note that $(\ref{conj9})$ gives
        \begin{flalign*}
        	&\lvert \rho{(\alpha,Ae^{f(\theta)})} - \frac{\la \deg B_{l_j}, \alpha \ra}{2} - \rho{(\alpha,A_{l_j})} \rvert \\
        	= & \lvert \rho{(\alpha,A_{l_j}e^{\overline f_{l_j}(\theta)})} - \rho{(\alpha,A_{l_j})} \rvert \\
        	\leqslant & \epsilon_{l_{j+1}}^{\frac{1}{2}}
        \end{flalign*}
        and $(\ref{estb5})$ implies 
        $$
        \deg B_{l_j} = \deg B_{l_{n_q}} = m_0.
        $$
        Therefore, we have
        \begin{equation}\label{estrot2}
        	\lvert \rho{(\alpha,A_{l_{n_q}}e^{\overline f_{l_{n_q}}(\theta)})} - \rho{(\alpha,A_{l_j})} \rvert \leqslant \epsilon_{l_{j+1}}^{\frac{1}{2}}.
        \end{equation}
        
        By $(\ref{estrot2})$ and $\rho(\alpha,A_{l_{n_q}}e^{\overline f_{l_{n_q}}(\theta)}) = 0$, for $m \in \Z^d$ with $0<|m| \leqslant N_{l_{j+1}}=\frac{2}{\frac{1}{l_{j+1}}-\frac{1}{l_{j+2}}} \lvert \ln \epsilon_{l_{j+1}}\rvert$, we have 
        \begin{flalign*}
        	&\lVert 2\rho(\alpha,A_{l_j})-\la m,\alpha \ra \rVert_{\R / \Z} \\
        	\geqslant &\lVert 2\rho(\alpha,A_{l_{n_q}}e^{\overline f_{l_{n_q}}(\theta)})-\la m,\alpha \ra \rVert_{\R / \Z}-| 2\rho(\alpha,A_{l_{n_q}}e^{\overline f_{l_{n_q}}(\theta)})-2\rho(\alpha,A_{l_j})| \\
        	\geqslant &\lVert \la m,\alpha \ra \rVert_{\R / \Z}-2\epsilon_{l_j}^\frac{1+s}{2} \\
        	\geqslant &\frac{\kappa}{|m|^ \tau} - 2\epsilon_{l_{j+1}}^\frac{1}{2} \\
        	\geqslant & \frac{\kappa}{\left( \frac{2}{\frac{1}{l_{j+1}}-\frac{1}{l_{j+2}}} \ln \frac{1}{\epsilon_{l_{j+1}}}\right) ^\tau} - 2\epsilon_{l_{j+1}}^\frac{1}{2} \\
        	\geqslant & 2\epsilon_{l_{j+1}}^\frac{\sigma}{2} - 2\epsilon_{l_{j+1}}^\frac{1}{2} \\
        	> &\epsilon_{l_{j+1}}^\sigma, \ \forall n_q+1 \leqslant j \leqslant j_0.
        \end{flalign*}
        This means the $(j_0 + 1)$-th step is also non-resonant with estimates
        $$
        B_{l_{j_0 + 1}}=\widetilde{B}_{l_{j_0}} \circ B_{l_{j_0}}, \ \lvert \widetilde B_{l_{j_0}}(\theta)- Id\rvert_{\frac{1}{l_{j_0 + 2}}} \leqslant \epsilon_{l_{j_0 + 1}}^{1-\frac{8}{D}}, \ \deg B_{l_{j_0}}= m_0.
        $$
        
        To conclude, for $\forall j > n_q$, we have 
        $$
        B_{l_j}(\theta+\alpha)(Ae^{f_{l_{j}}(\theta)})B^{-1}_{l_j}(\theta)=A_{l_j}e^{f_{l_j}^{'}(\theta)},
        $$
        with estimates
        \begin{equation}\label{estb6}
        	B_{l_{j}}=\widetilde{B}_{l_{j-1}} \circ B_{l_{j-1}}, \ \lvert \widetilde B_{l_{j-1}}(\theta)- Id\rvert_{\frac{1}{l_{j+1}}} \leqslant \epsilon_{l_j}^{1-\frac{8}{D}}, \ \deg B_{l_{j-1}} = m_0.
        \end{equation}
        
        Denote $B_2=\lim \limits_{j\to\infty}B_{l_j}$, $\tilde A_2=\lim \limits_{j\to\infty}A_{l_j} \in SL(2,\R)$. The rest of the process is similar to the above Diophantine case and then $(\ref{trans2})$ and $(\ref{estrot})$ holds. 
        
        This finishes our proof.
	\end{pf}
	
	\begin{Remark}
		Compared to Cai and Ge \cite{cai2022reducibility}, the reducibility results in this paper are certainly optimal. We have conducted a process parallel to the almost reducibility procedure. In contrast to the previous work, the conclusions in this paper are sharper in terms of regularity loss and the loss of regularity is independent of the parameter $k$.
	\end{Remark}

	\section{Spectral application}
	With the above reducibility theorem in hand, we can prove our spectral applications Theorem \ref{main2} and Theorem \ref{main3} as stated in the introduction from this. Let us first introduce several definitions and cite some results shown by \cite{ge2023ballistic}\cite{cai2021polynomial}.
	
	For spatial transport properties of a quantum particle on the lattice $\Z$, we are interested in studying the observable quantity associated with the position of the particle, which is represented by the unbounded self-adjoint operator 
	$$
	(Xx)_n:=(nx)_n, \ n \in \Z, 
	$$
	with its natural domain of definition 
	$$
	{\rm Dom} X = \left\{ x_n \in \ell^2(\Z):\sum\limits_{n \in \Z} |n|^2|x_n|^2 < +\infty \right\}. 
	$$
	
	Our focus is on the phenomenon of ballistic motion, which informally states that the particle's position grows linearly with time ($X(T) \approx T$). More precisely, we aim to investigate the following limit: 
	\begin{equation}\label{dio3}
		\lim\limits_{T \to +\infty} \frac{1}{T} X(T) x, 
	\end{equation}
	where the initial state at time $x \in {\rm Dom} X$. The limit (\ref{dio3}) can be regarded as the ``asymptotic velocity" of the state $x$ as time approaches infinity, provided that the limit exists. Then we can define the asymptotic velocity operator as 
	\begin{equation}\label{dio4}
		Q=\lim\limits_{T \to +\infty} \frac{1}{T} X(T) = \lim\limits_{T \to +\infty} \frac{1}{T} { \displaystyle \int_{0}^{T} e^{itH}Se^{-itH} dt}, 
	\end{equation} 
	where $S$ is bounded and satisfies
	\begin{equation}\label{dio5}
		Sx_n = i(x_{n+1} - x_{n-1}).  
	\end{equation}
	
	The Schrödinger operator $H$ is said to demonstrate strong ballistic transport if the strong limit on the right-hand side of equation (\ref{dio4}) exists, is defined on the entire $\ell^2(\Z)$ and $\ker Q = \{0\}$.
	
	\begin{Definition}\cite{ge2023ballistic}
	Let $\mathcal{K} \subset \R$ be a Borel subset. We call $H$ has strong ballistic transport on $\mathcal{K}$ if there exists a self-adjoint operator $Q$ such that 
	    \begin{equation}\label{dio6}
		    \lim\limits_{T \to +\infty} \frac{1}{T} { \displaystyle \int_{0}^{T} e^{itH}1_{\mathcal{K}}(H)S1_{\mathcal{K}}(H)e^{-itH} dt} = 1_{\mathcal{K}} Q 1_{\mathcal{K}}   
	    \end{equation}
	and	$\ker Q = {\rm Ran} (1_{\mathcal{K}})^{\perp}$, where $1_{\mathcal{K}}(\cdot)$ denotes the indicator function of $\mathcal{K}$.    
	\end{Definition}
	
	\begin{Definition}\cite{ge2023ballistic}
		Cocycle $(\alpha, S^V_E)$ is said to be $C^k$-reducible in expectation on $\mathcal{K}$ if it is $C^k$-reducible for each $E \in \mathcal{K}$ and there exists a choice of $L^2(\T^d)$-normalized conjugations $B(E; \cdot)$ such that 
		$$
		 \displaystyle \int_{\mathcal{K}} ||B(E; \cdot)||^4_{C^k(\T^d)}d\rho(E) < +\infty. 
		$$ 
	\end{Definition}
	
	\begin{Lemma}\cite{ge2023ballistic}
		Let a (Borel) subset $\mathcal{K} \subset \R$, $\{H_{V,\alpha,\theta} \}_{\theta \in \T^d}$ be a quasiperiodic operator family whose cocycles are $C^k$-reducible in expectation on $\mathcal{K}$ for some $k > \frac{5d}{2}$. Then the family $\{H_{V,\alpha,\theta} \}_{\theta \in \T^d}$ has strong ballistic transport on $\mathcal{K}$.
	\end{Lemma}
	
	As a direct application of our conclusion, now we can prove Theorem \ref{main2} in the following. 
	
	\begin{pf}    
		We will not consider $E$ in the spectral gap, because cocycle $(\alpha, S^V_E)$ is uniformly hyperbolic in the quasi-periodic Schr\"odinger case. Therefore, it is always reducible under our conditions and the conclusion is naturally valid. 
		
		For $E \in \Sigma_{V,\alpha}$, Theorem \ref{thm3.1} shows that cocyle $(\alpha, S^V_E)$ is reducible. As $k > 14\tau + 2$ and picking $k_0 = [k-10\tau-3]$, we can get $k_0 > [4\tau-1]$ and $\tau>d$, then condition $k_0 > \frac{5d}{2}$ is obviously true. Since $E$, which satisfies neither the Diophantine nor rational condition in the spectrum set, is a set of measure zero. Note that, in Ge and Kachkovskiy \cite{ge2023ballistic}, under the setting of Zhao \cite{zhao2016ballistic}, we have strong ballistic transport on the whole spectrum. In the case of $C^k$ in this paper, the conclusion still holds. Therefore, by Lemma 4.1, the Schr\"{o}dinger operator $H_{V,\alpha,\theta}$ has strong ballistic transport for a.e. $E \in \Sigma_{V,\alpha}$. 
		
		Readers can also refer to the detailed proof in Ge and Kachkovskiy \cite{ge2023ballistic}. 
		
		This finishes the proof of Theorem 1.2. 
	\end{pf}
	
	Now, for the application of the reducibility theorem to spectral structures Theorem \ref{main3}, we give the proof as follows.  
	
	\begin{pf}
	    Referring to Cai and Wang \cite{cai2021polynomial}, we know the homogeneity of the spectrum is connected with polynomial decay of gap length and H\"{o}lder continuity of IDS. Then we need the following lemmas. 
	    
	    \begin{Lemma}\cite{cai2021polynomial}
	    	Let $G_m(V)=(E^-_m, E^+_m)$ denote the gap with label $m$ and, $\alpha \in {\rm DC}_d(\kappa, \tau), V \in C^k(\T^d, \R)$ with $k \geqslant D_0 \tau$, where $D_0$ is a numerical constant. There exists $\tilde \epsilon = \tilde \epsilon (\kappa, \tau, k, d) > 0$ such that if $\lVert V \rVert_k \leqslant \tilde \epsilon$, then $\lvert G_m(V) \rvert_k \leqslant {\tilde \epsilon}^{\frac{1}{4}}\lvert m \rvert^{- \frac{k}{9}}$. 
	    \end{Lemma}
	    
	    \begin{Lemma}\cite{cai2022absolutely}
	    	Let $\alpha \in {\rm DC}_d(\kappa, \tau), V \in C^k(\T^d, \R)$ with $k>17\tau+2$, then there exists $\lambda_0$ depending on $V, d, \kappa, \tau, k$ such that if $\lambda < \lambda_0$, then $N_{\lambda V, \alpha}$ is $\frac{1}{2}$-H\"{o}lder continuous: 
	    	$$
	    	N(E+\hat \epsilon)-N(E-\hat \epsilon) \leqslant C_0 \hat \epsilon^{\frac{1}{2}}, \ \forall \hat \epsilon >0, \ \forall E \in \R. 
	    	$$
	    	where $C_0$ depends only on $d, \kappa, \tau, k$.  
	    \end{Lemma}
	    
	    Through the quantitative almost reducibility theory and reducibility theory, we can show both of the above lemmas. Under the preconditions of Theorem \ref{thm3.1} together with lemma 4.2 and lemma 4.3, for any $E \in \Sigma_{V,\alpha}$ with $2\rho(\alpha, S^V_E) =\la m_0,\alpha  \ra \mod \Z$ for some $m_0\in \Z^d$, we will divide $E$ into three cases when considering $\nu$-homogeneous. In each case, by calculating directly, the Definition 1.1 is satisfied for every $E \in \Sigma_{V,\alpha}$. Please refer to the detailed proof in Cai and Wang \cite{cai2021polynomial}. 
	    
	    This finishes the proof of Theorem 1.3. 	 
	\end{pf}

	\section{Acknowledgements}
	Ao Cai was supported by NSFC grant 12301231. 
	
	Zhiguo Wang was supported by NNSF grant 12071410 and 12071232.

	\bibliographystyle{amsplain}
	\bibliography{ref}
	
\end{document}